\documentclass[a4paper,10pt]{amsart}

\usepackage{amsmath,enumerate,amsfonts,amssymb,amsthm}
\usepackage[dvips]{graphicx}

\theoremstyle{plain}
\newtheorem*{Theo}{\bf THEOREM}
\newtheorem{Prop}{\bf PROPOSITION}
\newtheorem{Lem}{\bf LEMMA}

\newtheorem{Def}{\bf DEFINITION}
\newtheorem{Ex}{\bf EXAMPLE}
\theoremstyle{remark}
\newtheorem*{Rem}{\bf Remark}

\DeclareMathOperator{\re}{Re} \DeclareMathOperator{\im}{Im}
 \DeclareMathOperator{\hess}{hess}
\DeclareMathOperator{\pf}{pf} \DeclareMathOperator{\tr}{trace}
 
\DeclareMathOperator{\vol}{vol} 
\DeclareMathOperator{\signature}{signature}

\begin{document}

\title{Integrable geometries and Monge-Amp\`ere
equations}
\author{Bertrand Banos}
\email{bertrand.banos@wanadoo.fr}

\maketitle

\begin{abstract}
In this lecture delivered at the Integrable and Quantum Field
Theory at Peyresq sixth meeting, we review the Lychagin's
Monge-Amp\`ere operators theory and exhibit the link it
establishes between the classical problem of local equivalence for
non linear partial differential equations and the problem of
integrability of some geometrical structures.
\end{abstract}

\section*{Introduction}

A general approach to the study of non-linear Partial Differential
Equations, which goes back to Sophus Lie, is to see a $k$-order
equation on a $n$-dimensional manifold $N^n$ as a closed subset in
the manifold of $k$-jets $J^kN$. In particular, a second-order
differential equation lives in the space $J^2N$. Nevertheless, as
it was noticed  by Lychagin in his seminal paper "Contact geometry
and non-linear second-order differential equations" (1979), it is
possible to decrease one dimension and to work on the contact
space $J^1N$ for a large class of second order PDE's, containing
quasi linear PDE's and Monge-Amp\`ere equations.

 Moreover, for a large class of operators (those
which admit a symmetry), we can replace the $1$-jet space by the
cotangent space and contact geometry by symplectic geometry. This
study of differential operators becomes then the study of
differential forms in the presence of a symplectic form.

The aim of this review paper is to use this Lychagin's
correspondence to show that it is possible to reconstruct the
geometrical background starting from such a partial differential
equation and to interpret the integrability of this geometry in
terms of "integrability" of the equation.

In the first part, a brief review on the Monge-Amp\`ere operator
theory is given. Concepts of generalized solution and local
equivalence are presented in geometric terms.

In the second part, the geometry of differential forms is studied.
A unified approach is given in any dimensions and classification
results in dimensions $2$ and $3$ are presented then.

In the last part, a link between $2D$-Monge-Amp\`ere equations and
generalized complex geometry is described. Conservation laws and
Generating functions are presented as generalized complex objects.

Many thanks to the organizers Paul Baird, Fr\'ed\'eric H\'elein,
Joseph Kouneiher, Franz Pedit and Volodya Roubtsov for their kind
invitation !

\section{Monge-Amp\`ere operators theory}

Let $M^n$ be a $n$-dimensional manifold, $T^*N$ its cotangent
bundle and $\Omega$ the canonical symplectic structure on it.
Locally,
$$
\Omega=\sum_{i=1}^n dq_i\wedge dp_i,$$ with $(q_1,\ldots q_n)$
coordinates on $M$. We denote by $\Omega^*(T^*M)$ the space of
differential forms on the $2n$-dimensional manifold $T^*M$. For
example, $\Omega\in \Omega^2(T^*M)$.

\subsection{Monge-Amp\`ere operators}

\begin{Def}
Let $\omega\in \Omega^n(T^*M)$. The Monge-Amp\`ere operator
associated with $\omega$ is the differential operator
$$\Delta_\omega: C^\infty(M)\rightarrow \Omega^n(M)$$
defined by
$$
\Delta_\omega(f)=(df)^*(\omega),
$$
with $df:M\rightarrow T^*M$ the natural section defined by $f$.
\end{Def}

\begin{Ex}
Consider on $T^*\mathbb{R}^2$, the $2$-form $\omega=dq_ 1\wedge
dp_2 -dq_2\wedge dp_1$. We get
$$
\Delta_\omega(f)=\big(f_{q_1q_1} + f_{q_2q_2}\big)dq_1\wedge dq_2.
$$
\end{Ex}

\begin{Ex}
Consider on $T^*\mathbb{R}^3$, the $3$-form
$$
\omega=dp_1\wedge dq_2\wedge dq_3 + dq_1\wedge dp_2\wedge dq_3 +
dq_1\wedge dq_2\wedge dp_3 - dp_1\wedge dp_2\wedge dp_3.
$$
We get
$$
\Delta_\omega(f)=\big(\Delta f - \hess(f)\big) dq_1\wedge
dq_2\wedge dq_3.
$$
\end{Ex}

We obtain a large class of non linear partial differential
equations, characterized by their "determinant like" nonlinearity.
This class is called the class of symplectic Monge-Amp\`ere
equations (SMAE). The term symplectic means that the corresponding
form lives on the cotangent bundle. The whole class of
Monge-Amp\`ere equations is obtained with forms on the contact
manifold $J^1M$.

Using this correspondence between SMAE and differential forms, we
will see now how one can describe in geometric terms two classical
notions in the study of PDE's: the notion of generalized solution
and the notion of local equivalence.

\subsection{Generalized solutions}

\begin{Def}
A generalized solution of a SMAE $\Delta_\omega=0$ is a lagrangian
submanifold $L^n$ of the symplectic manifold $(T^*M,\Omega)$ on
which vanishes $\omega$:
$$
\omega|_L=0.
$$
\end{Def}

A lagrangian submanifold which is a graph is the graph of a closed
form. Hence, a generalized solution can be thought as a smooth
patching of regular solutions. \vspace{1cm}
\begin{figure}[hbp!]
\includegraphics{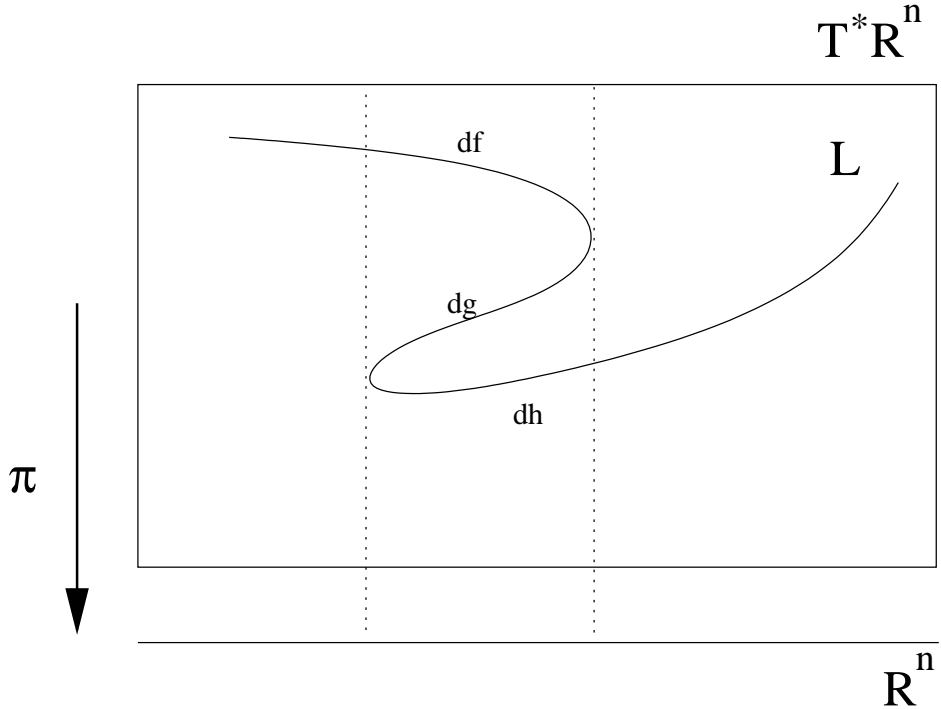}
\caption{Generalized solution}
\end{figure}

\begin{Ex}
On $\mathbb{R}^2$, a regular solution of the Laplace equation is
of course an harmonic function. A generalized solution is a
surface of $\mathbb{C}^2$ on which vanish $\Omega=\re(dz_1\wedge
dz_2)$ and $\omega=\im(dz_1\wedge dz_2)$, that is a complex curve
of $\mathbb{C}^2$.
\end{Ex}

\begin{Ex}
On $\mathbb{R}^3$, a generalized solution of the SMAE
$$
\Delta f - \hess(f)=0,
$$
is a submanifold of $T^*\mathbb{R}^3=\mathbb{C}^3$ on which vanish
the symplectic form
$$
\Omega=\frac{i}{2}\sum_{j=1}^3 dz_j\wedge d\bar{z_j}
$$
and the $3$-form
$$
\omega=\im(dz_1\wedge dz_2\wedge dz_3),
$$
that is, a special lagrangian submanifold of $\mathbb{C}^3$.
\end{Ex}

For any form $\theta\in \Omega^{n-2}(T^*M)$, the form
$\theta\wedge\Omega$ vanished on all lagrangian submanifolds of
$T^*M$. The correspondence between SMAE and differential forms is
therefore well defined up these particular forms. We introduce
then the notion of primitive form:

\begin{Def}
A $n$-form $\omega\in \Omega^n(T^*M)$ is said to be primitive if
$$
\omega\wedge \Omega=0.
$$
\end{Def}

\begin{Theo}[Hodge-Lepage-Lychagin]
\begin{enumerate}
\item any differential form $\omega\in \Omega^n(T^*M)$ admits an
unique decomposition
$$
\omega=\omega_0+ \omega_1\wedge\Omega
$$
with $\omega_0$ primitive.
\item if two primitive $n$-forms vanish on the same lagrangian
submanifolds, then they are proportional.
\end{enumerate}
\end{Theo}

\subsection{Local equivalence}

Roughly speaking, two PDE's are locally equivalent if they have
the same solutions, up a change of dependent and independent
coordinates. In Lychagin's formalism, this notion becomes clear:

\begin{Def}
Two SMAE $\Delta_{\omega_1}=0$ and $\Delta_{\omega_2}=0$ are said
to be (locally) equi\-valent if there exists a (local)
diffeomorphism $F:T^*M\rightarrow T^*M$ preserving the symplectic
form, that is
$$
F^*\Omega=\Omega,
$$
and exchanging the two forms, that is
$$
F^*\omega_1=\omega_2.
$$
\end{Def}

Note that $F$ sends a generalized solution of
$\Delta_{\omega_2}=0$ on a generalized solution of
$\Delta_{\omega_1}=0$, but not necessarily a regular solution on a
regular solution.

\begin{Ex}
Let us consider on $\mathbb{R}^2$, the Monge-Amp\`ere equation
$$
f_{q_1q_1}f_{q_2q_2}-f_{q_1q_2}^2=1,
$$
corresponding to the form
$$
\omega=dp_1\wedge dp_2-dq_1\wedge dq_2.
$$
Let $\phi:T^*\mathbb{R}^2\rightarrow T^*\mathbb{R}^2$ be the
partial Legendre transform
$$
\phi(q_1,q_2,p_1,p_2)=(q_1,p_2,p_1,-q_2).
$$
Since
$$
\phi^*(\omega)=dq_2\wedge dp_1 - dq_1\wedge dp_2,
$$
we see that our Monge-Amp\`ere equation is equivalent to the
Laplace equation. In theory, we can construct solutions using
harmonic functions on $\mathbb{C}$. For example, consider the
harmonic function $f(q_1,q_2)=e^{q_{1}}\sin(q_2)$. One can check
that the graph $L_f$ of $df$ is sent by $\phi$ on a submanifold
which is itself a graph:
$$
\phi(L_f)=\big\{t_1,t_2,g_{t_1},g_{t_2}\big)\}
$$
with
$$
g(t_1,t_2)=q_2\arcsin(q_2e^{-q_1})+ \sqrt{e^{2q_1}-q_2^2}.
$$
This function $g$ is a non trivial solution of our Monge-Amp\`ere
equation.
\end{Ex}

\section{Geometry of differential forms}

Hence, the classical problem of local equivalence for
Monge-Amp\`ere equations can be understood  as a problem of the
Geometric Invariant Theory: the idea is to construct invariant
structures which will characterize each equivalent class.

The first step of this approach is pointwise: we study the action
of the symplectic group $SP(n\mathbb{R})$ on the space of
primitive forms $\Lambda^n_0(\mathbb{R}^n)$. We will see next how
one can "integrate" such study.

\subsection{The bracket}

Let $V^{2n}$ be a $2n$ dimensional  real vector space. We fix a
symplectic form $\Omega$ on $V$ and the volume form
$$\vol=\frac{\Omega^n}{n!}.$$

We denote by $\Lambda^n(V^*)$ the space of $n$-forms on $V$ and by
$\Lambda^n_0(V^*)$ the space of primitive $n$-forms, that is
$$
\Lambda_0^n(V^*)=\{\omega\in \Lambda^n(V^*),\;\Omega\wedge\omega=0\}.
$$

We denote by $SL(2n)$ the group of automorphisms preserving the
volume form $\vol$ and by $SP(n,\mathbb{R})$ the group of
automorphisms preserving the symplectic form $\Omega$. Their Lie
algebras are denoted by $sl(2n)$ and $sp(n,\mathbb{R})$.

Using the exterior product, we define an isomorphism $\mathcal{A}:
\Lambda^{2n-1}(V^*)\rightarrow V$ by
$$
<\alpha,
\mathcal{A}(\theta)>=\frac{\alpha\wedge\theta}{\vol},\;\;\;
\text{for $\alpha \in \Lambda^1(V^*)$ and $\theta\in
\Lambda^{2n-1}(V^*)$.}
$$

\begin{Def}
The bracket $\Phi:\Lambda^n(V^*)\times \Lambda^n(V^*)\rightarrow
sl(V)$ is defined by
$$
\Phi(\omega_1,\omega_2)(X)=\mathcal{A}\Big((\iota_X\omega_1)\wedge
\omega_2 - (-1)^n \omega_1\wedge (\iota_X\omega_2)\Big).
$$
\end{Def}

It is straightforward to check the two following lemmas:

\begin{Lem}
This bracket is invariant under the action of $SL(2n)$, that is
$$
\Phi(F^*\omega_1,F^*\omega_2)=F^{-1}\circ
\Phi(\omega_1,\omega_2)\circ F
$$
for any $F\in SL(2n)$.
\end{Lem}

\begin{Lem}{\label{rec}}
Let $\tilde{\Phi}$ be the bracket defined for the $(2n+2)$
dimensional vector space $\tilde{V}= V\times
\mathbb{R}_{t_1}\times \mathbb{R}_{t_2}$ endowed with the volume
form
$$
\tilde{\vol}=\vol\wedge dt_1\wedge dt_2.$$ Then the following
relations hold
\begin{enumerate}
\item $\tilde{\Phi}(\omega_1\wedge dt_1,\omega_2\wedge
dt_2)(\partial_{t_1})=-\tilde{\Phi}(\omega_1\wedge
dt_1,\omega_2\wedge dt_2)(\partial_{t_2})$
\item $\tilde{\Phi}(\omega_1\wedge dt_1,\omega_2\wedge
dt_2)(X) = \Phi(\omega_1,\omega_2)(X),\;\;\; \forall X\in V.$
\end{enumerate}
\end{Lem}

Note that this second lemma shows that $\Phi$ takes  its values in
$sl(2n)$.

In the case $n=3$, the tensor
$K_\omega=\frac{1}{2}\Phi(\omega,\omega)$ is the invariant
constructed by Hitchin in \cite{H1}, which can be easily extended
to any odd $n$:

\begin{Prop}[Hitchin]{\label{moment}} When $n$ is odd, the map
$K:\Lambda^n(V^*)\rightarrow sl(2n)$, $\omega\mapsto
\frac{1}{2}\Phi(\omega,\omega)$ is a moment map for the
hamiltonian action of $SL(2n)$ on $\Lambda^n(V^*)$ endowed with
the symplectic form
$$
\Theta(\omega_1,\omega_2)=\frac{\omega_1\wedge\omega_2}{\vol}.
$$
\end{Prop}

When $n$ is even, the bracket $\Phi$ is antisymmetric and the
situation is completely different. The analog of \ref{moment} is
the following, which is proved in \cite{BR} (in preparation):

\begin{Prop}{\label{Lie}}
We define on $\Lambda^n(V^*)\times sl(2n)$ the following bracket:
\begin{enumerate}
\item $[A_1,A_2]=A_1A_2-A_2A_1$
\item $[A,\omega]=L_A(\omega)$
\item $[\omega_1,\omega_2]=\Phi(\omega_1,\omega_2)$
\end{enumerate}
for $A$, $A_1$ and $A_2$ in $sl(2n)$ and $\omega$, $\omega_1$ and
$\omega_2$ in $\Lambda^n(V^*)$.

Then $[\;,\;]$ is a Lie bracket.
\end{Prop}

The $SP(n,\mathbb{R})$-version of these results is summed up in
the following:

\begin{Prop}
\begin{enumerate}
\item if $\omega_1$ and $\omega_2$ are primitive then
$\Phi(\omega_1,\omega_2)\in sp(n,\mathbb{R})$.
\item if $n$ is odd, then $K:\Lambda^n_0(V^*)\rightarrow
sp(n,\mathbb{R})$, $\omega\mapsto \frac{1}{2}\Phi(\omega,\omega)$
is a moment map for the hamiltonian action of $SP(n,\mathbb{R})$
on the symplectic subspace $\Lambda^n_0(V^*)$ of $\Lambda^n(V^*)$.
\item if $n$ is even, the space $\Lambda^n_0(V^*)\oplus
sp(n,\mathbb{R})$ is a Lie subalgebra of $\Lambda^n(V^*)\oplus
sl(2n)$.
\end{enumerate}
\end{Prop}

\begin{Rem}
When $n$ is odd, the tensor $K_\omega$ defines a family of scalar
invariants
$$a_k=\tr(K_\omega^{2k}),\;\;\; k\in\mathbb{N}$$
and a quadratic form called the Lychagin-Roubtsov quadratic form:
$$
q_\omega(X)=\Omega(K_\omega X,X).
$$

When $n$ is even, the adjoint operator $ad_\omega=[\omega,\cdot]$
defines an endomorphism
$$ad_\omega^2: sp(n,\mathbb{R})\rightarrow
sp(n,\mathbb{R})$$ which gives also a family of scalar invariants
$$
a_k=\tr(ad_\omega^{2k}),\;\;\; k\in\mathbb{N}
$$
and a symmetric polynomial of degree $4$ defined by
$$
q_\omega(X)=\tr(\lbrack
ad_\omega^2(X\otimes\iota_X(\Omega))\rbrack^2).
$$
\end{Rem}

\subsection{Examples}

\subsubsection{$n=2$}

The identity  $\omega= \Omega(A_\omega\cdot,\cdot)$ gives an
isomorphism between the space of $2$-forms
$\Lambda^2(\mathbb{R}^4)$ and the Jordan algebra $Jor(\Omega)$
defined by
$$
Jor(\Omega)=\{A\in gl(4),\;
\Omega(A\cdot,\cdot)=\Omega(\cdot,A\cdot)\}.
$$
Our bracket $\Phi$ becomes then the usual bracket:
$$
\Phi(\omega_1,\omega_2)=A_{\omega_1}A_{\omega_2}-
A_{\omega_{2}}A_{\omega_1}.
$$
We easily see then the isomorphism of Lie algebras
$$
\Lambda^2_0(\mathbb{R}^4)\oplus sp(2,\mathbb{R}) =
sl(4,\mathbb{R}).
$$

Moreover, for $\omega\in \Lambda_0^2(\mathbb{R}^4)$, the
endomorphism $ad^2_\omega=ad^2_{A_\omega}:
sp(2,\mathbb{R})\rightarrow sp(2,\mathbb{R})$ satisfies
$$
\tr(ad^2_\omega) = 16\pf(\omega)
$$
where the pfaffian of $\omega$ is the classical invariant
$$
\pf(\omega)=\frac{\omega\wedge\omega}{\Omega\wedge\Omega}.
$$
The polynomial $q_\omega$ is the null polynomial.

\subsubsection{$n=3$.}

It is proved in \cite{H1} that the action of $GL(6,\mathbb{R})$ on
$\Lambda^3(\mathbb{R}^3)$ has two opened orbits separated by the
hypersurface $\lambda=0$ where
$$
\lambda(\omega)=\frac{1}{6}\tr(K_\omega^2).
$$
Note that, for any $3$-form the following holds:
$$
K_\omega^2=\lambda(\omega)\cdot Id.
$$

By analogy with the $2$-dimensional case, we call this invariant
the Hitchin pfaffian. A $3$-form with a non vanishing Hitchin
pfaffian is said to be nondegenerate.

For a primitive form $\omega$ , we get a triple $(g_\omega,
K_\omega,\Omega)$ with $g_\omega=\Omega(K_\omega\cdot,\cdot)$ the
Lychagin-Roubtsov metric (see \cite{LR}, \cite{Ba2}).  This triple
defines  a  $\epsilon$-K\"ahler structure in the sense of
\cite{S},  that is,  the tensor $K_\omega$ satisfies, up a
renormalization,
$$
K_\omega^2=\epsilon,\text{ with $\epsilon =0,1,-1$}.
$$
Note that the Lychagin-Roubtsov metric has signature.

Moreover, in the nondegenerate case the form $\omega$ admits an
unique dual form $\hat{\omega}$, such that
$\omega+\sqrt{\epsilon}\hat{\omega}$ and
$\omega-\sqrt{\epsilon}\hat{\omega}$ are the volume forms of the
two-eigenspaces of the Hitchin tensor $K_\omega$. Saying
differently, to each nondegenerate primitive forms corresponds a
$\epsilon$-Calabi-Yau structure.

\subsubsection{$n=4$}

The Lie algebras $\Lambda^4(\mathbb{R}^8)\oplus sl(8,\mathbb{R})$
and $\Lambda_0^4(\mathbb{R}^8)\oplus sp(4,\mathbb{R})$ are known
to be isomorphic to the exceptional Lie algebras $E_7$ and $E_6$
(see \cite{W}). Moreover, it is proved in \cite{K} that the family
$\{a_k=\tr(ad_\omega^{2k})\}_{k\in\mathbb{N}}$ forms a complete
family of invariants.

Nevertheless, computations in these dimensions are extremely
complicated. The author planes to implement an algorithm which
could give in a reasonable time these invariants $a_k$.

It is worth mentioning that, on many examples, the symmetric
polynomial $q_\omega$ of degree $4$  is the square of a quadratic
form. Is it always true ? A positive answer would be extremely
useful to understand the geometry of PDE's of Monge-Amp\`ere type
in $4$ variables.

\subsection{Classifications results}

\subsubsection{Monge-Amp\`ere equations in $2$ and $3$ variables}

The action of the symplectic linear group on $2D$ and $3D$
symplectic Monge-Amp\`ere equations with constant coefficients has
a finite numbers of orbits and we know all of them as it is shown
in tables 1 and 2 (see \cite{LR} and \cite{Ba1}).

\begin{table}[ht!]
\begin{center}
\begin{tabular}{|c|c|c|}
\hline \mathversion{bold} $\Delta_\omega=0$&\mathversion{bold}
$\omega$ & \mathversion{bold} $\pf(\omega)$\\ \hline \hline
$\Delta f=0$ & $dq_1\wedge dp_2 - dq_2\wedge dp_1$&$1$\\
\hline
$\square f=0$& $dq_1\wedge dp_2 + dq_2\wedge dp_1$&$-1$\\
\hline $\frac{\partial^2 f}{\partial q_1^2}=0$ & $dq_1\wedge dp_2$
& $0$\\
 \hline
\end{tabular}
\caption{Classification of SMAE in $2$ variables}
\end{center}
\end{table}

\begin{table}[!ht]
\begin{center}
\begin{tabular}{|c|c|c|c|}
\hline &\mathversion{bold}$\Delta_\omega=0$ &
\mathversion{bold}$\signature(q_\omega)$ &
\mathversion{bold}$\lambda(\omega)$\\ \hline \hline
 1& $\hess(f)=1$ &$(3,3)$&$1$\\
\hline
 2& $\Delta f- \hess(f)=0$ &$(0,6)$&$-1$\\
\hline 3& $\square f +\hess(f)=0$ &$(4,2)$&$-1$\\
\hline \hline
 4& $\Delta f=0$&$(0,3)$&$0$\\
\hline
5&$\square f=0$&$(2,1)$&$0$\\
\hline
6&$\Delta_{q_2,q_3} f=0$&$(0,1)$&$0$\\
\hline
7&$\square_{q_2,q_3} f=0$&$(1,0)$&$0$\\
\hline
8&$\frac{\partial^2 f}{\partial q_1^2}=0$&$(0,0)$&$0$\\
\hline
\end{tabular}
\caption{Classification of SMAE in $3$ variables}
\end{center}
\end{table}

\begin{Rem}
\begin{enumerate}
\item in two variables, any SMAE with constant coefficients is
linearizable, that is equivalent to a linear PDE. Moreover, the
pfaffian distinguishes the different orbits.
\item in three variables, there exist nonlinearizable SMAE with
constant
coefficients and they correspond to nondegenerate primitive
$3$-forms. Moreover, the Hitchin pfaffian does not distinguish the
different orbits but so does the signature of the
Lychagin-Roubtsov metric.
\end{enumerate}
\end{Rem}

In $4$ variables, the action of the symplectic group is not
discrete anymore and there is no hope to obtain an exhaustive list
as in $2$ or $3$ variables. Moreover, it appears that on many
interesting examples, the associated geometry is completely
degenerated. Saying differently, in $4$ variables appears the
notion of non linear but degenerated Monge-Amp\`ere equation. In
the table 3, we have computed the polynomial invariant $q_\omega$
for the follwing examples:
$$
\hess(u)=1 \;\;\;\;\;\text{(Usual Monge-Amp\`ere equation)}
$$
$$
 \hess(u) - (\sum_{i<j}
u_{q_iq_i}u_{q_jq_j}-u_{q_iq_j}^2) + 1 =0 \;\;\;\;\;\text{(4D SLAG
equation)}
$$
$$
u_{q_1q_2}u_{q_3q_4}-u_{q_1q_4}u_{q_2}q_{3}=1\;\;\;\;\;
\text{(Plebanski I equation)}
$$
$$
u_{q_1q_1}u_{q_3q_3}-u_{q_1q_3}^2+u_{q_1q_2}-u_{q_3q_4}=0
\;\;\;\;\; \text{(Plebanski II equation)}
$$
$$
u_{q_1q_1}+u_{q_1q_4}u_{q_2q_3}-u_{q_1q_3}u_{q_2q_4}=0
\;\;\;\;\;\text{(Grant equation)}
$$

\begin{table}[ht!]
\begin{center}
\begin{tabular}{|c|c|}
\hline \mathversion{bold} $\Delta_\omega=0$& \mathversion{bold}
$q_\omega$\\ \hline \hline usual Monge-Amp\`ere & $(dq_1dp_1 +
dq_2dp_2+dq_3dp_3+ dq_4dp_4)^2$\\
\hline SLAG&
$(dq1^2+dq_2^2+dq_3^2+dq_4^2+dp_1^2+dp_2^2+dp_3^2+dp_4^2)^2$\\
\hline Plebanski I & $0$ \\
\hline Plebanski II & $dq_1^4$\\
\hline Grant  & $0$\\
\hline
\end{tabular}
\caption{Examples of SMAE in $4$ variables}
\end{center}
\end{table}

\subsubsection{Note on ellipticity of Monge-Amp\`ere equations}

Recall that a second order linear partial differential equation
$$
\sum_{i,j=1}^n A_{ij}\frac{\partial^2 u}{\partial x_i\partial
x_j}=0
$$
is said to be
\begin{enumerate}
\item elliptic if the symmetric matrix $A$ has signature $(n,0)$
or $(0,n)$,
\item hyperbolic if the symmetric matrix $A$ has signature $(n-1,1)$
or $(1,n-1)$,
\item parabolic if the symmetric matrix $A$ is degenerate.
\end{enumerate}

Following Harvey and Lawson (\cite{HL})), we will say that a
Monge-Amp\`ere equation is elliptic, hyperbolic or parabolic if it
is so in a first order approximation.  More precisely, let
$\Delta_\omega=0$ be a Monge-Amp\`ere equation and $\phi$ a
solution. The linearization of $\Delta_\omega$ at $\phi$ is the
linear differential operator
$$
D_\phi(\Delta_\omega)(u)=\frac{d}{dt}\Big{|}_{t=0}
\Delta_\omega(\phi+tu).
$$
The equation $\Delta_\omega=0$ is said to be elliptic, hyperbolic
or  parabolic at the point $\phi$ if its linearization
$D_\phi(\Delta_\omega)=0$ is elliptic, hyperbolic or  parabolic.

\begin{Ex}
Let us consider a generic 2D SMAE
$$
A + B\psi_{x x } + 2C\psi_{x y } + D\psi_{y y } + E(\psi_{x x
}\psi_{y y } - \psi_{x y }^2) =  0.
$$
Its linearization at $\phi$ is
$$
(B+E\phi_{yy})u_{xx}+
2(C-E\phi_{xy})u_{xy}+(D+E\phi_{xx})u_{yy}=0.
$$
and since
$$
\begin{aligned}
&
\begin{vmatrix}
B+E\phi_{yy}& C-E\phi_{xy}\\
C-E\phi_{xy}&D+E\phi_{xx}\\
\end{vmatrix}\\
&= BD-C^2 +E\big(B\phi_{x x } + 2C\phi_{x y } + D\phi_{y y } +
E(\phi_{x x
}\phi_{yy } - \phi_{x y }^2)\big)\\
& =  BD-C^2-AE,
\end{aligned}
$$
we deduce that our $2D$-SMAE is elliptic for instance if and only
$BD-C^2-AE>0$, that is if and only if its corresponding primitive
form has positive pfaffian everywhere.
\end{Ex}

\begin{Ex}
Let us consider now the 3D special lagrangian equation
$$
\Delta \psi-\det(\text{\emph{Hess}}\;\psi)=0.
$$
Its linearization at a point $\phi$ is
$\overset{3}{\underset{i,j=1}{\sum}} A_{ij}\frac{\partial^2
u}{\partial x_i\partial x_j}=0$ with
$$
A=I_3-\Phi^*
$$
where $I_3$ denotes the matrix identity, $\Phi$ denotes the
hessian matrix of $\phi$ and $\Phi^*$ denotes its comatrix (the
matrix of cofactors).

Choose a basis in which $\Phi=\begin{pmatrix}
 \lambda_1 & 0 & 0\\
 0 & \lambda_2 & 0\\
 0 & 0 & \lambda_3\\
\end{pmatrix}$ with
$\lambda_1+\lambda_2+\lambda_3=\lambda_1\lambda_2\lambda_3$. You
get
$$
\begin{aligned}
A&= \begin{pmatrix}
1-\lambda_2\lambda_3 & 0 & 0\\
 0 & 1-\lambda_1\lambda_3 & 0\\
 0 & 0 & 1-\lambda_1\lambda_2\\
\end{pmatrix}\\
&=
(1-\lambda_1\lambda_2-\lambda_1\lambda_3-\lambda_2\lambda_3)\begin{pmatrix}
\frac{1}{1+\lambda_1^2}&0&0\\0&\frac{1}{1+\lambda_2^2}
&0\\0&0&\frac{1}{1+\lambda_3^2}\end{pmatrix}
\end{aligned}
$$
The 3D special lagrangian equation is therefore elliptic
everywhere.
\end{Ex}

The following formula generalizes this last example and gives a
relation between the linearization of a $3D$-SMAE, the
Lychagin-Roubtsov metric and, surprisingly, its dual equation:

\begin{Prop}The linearisation at a point $\phi$ of a $3D$-SMAE
$\Delta_\omega=0$ is
$$D_\phi(\Delta_\omega)(u)=\sum_{i,j=1}^3
B_{ij}\frac{\partial^2 u}{\partial x_i\partial x_j}
$$
with
$$
B=\Delta_{\hat{\omega}}(\phi)\cdot\big(g_\omega^{-1}\big)|_{L_\phi}.
$$
\end{Prop}

\subsection{Integration of the classification}

After this linear approach,the next step is to try to integrate
the classification: when is a given MAE $\Delta_\omega=0$
equivalent to a MAE with constant coefficient
$\Delta_{\omega_c}=0$ ? One can try to understand this classical
problem of integrability in terms of integrability of a certain
geometric structure. The idea is that a MAE contains information
on its underlying geometry. For example, $2D$-Laplace equation
contains information about the complex structure of $\mathbb{R}^4$
and $3D$-special lagrangian equation contains information about
the Calabi-Yau structure of $\mathbb{C}^3$. Using the
Lychagin-Roubtsov metric $g_\omega$ and the Hitchin tensor
$K_\omega$, one can actually define for any symplectic $3D$
Monge-Amp\`ere equation some geometrical structure of Calabi-Yau
type. The following result is proved in \cite{Ba2}:

\begin{Theo}
A symplectic $3D$ MAE is locally equivalent to one of the
following equations
$$
\begin{aligned}
&\hess(f)=1\\
&\Delta f - \hess(f)=0\\
&\square f + \hess(f)=0\\
\end{aligned}
$$
if and only if the structure of Calabi-Yau type it defines is
nondegenerate, flat and integrable.
\end{Theo}

This result has to be compared with its $2D$ analog obtained in
\cite{LR}:

\begin{Theo}
A symplectic $2D$ MAE is locally equivalent to one of the
following equation
$$
\begin{aligned}
&\Delta f =0\\
&\square f=0\\
\end{aligned}
$$
if and only the almost complex structure or almost product
structure it defines is integrable.
\end{Theo}

\section{2D Monge-Amp\`ere equations of divergent type and generalized
complex
geometry}

These results are quite frustrating: which kind  of integrable
geometries could we define for more general MAE ? One answer could
be: generalized complex geometry.

This very rich concept defined recently by Hitchin (\cite{H2}) and
developed by Gualtieri (\cite{Gu1}), which interpolates between
complex and symplectic geometry, is very popular since it seems to
provide a well-adapted geometric framework for different models in
string theory.

\subsection{Monge-Amp\`ere equations of divergent type}

Let us introduce first the Euler operator and the notion of
Monge-Amp\`ere equation of divergent type (see \cite{L}).
\begin{Def}
The Euler operator is the second order differential operator
$\mathcal{E}: \Omega^2(M)\rightarrow \Omega^2(M)$ defined by
$$
\mathcal{E}(\omega)=d \bot d\omega.
$$
A Monge-Amp\`ere equation $\Delta_\omega=0$ is said to be of
divergent type if $\mathcal{E}(\omega)=0$.
\end{Def}

\begin{Ex}[Born-Infeld Equation]
The Born-Infeld equation is
$$
(1-f_t)^2 f_{xx}+2f_tf_xf_{tx} - (1+f_x^2)f_{tt}=0.
$$
The corresponding primitive form is
$$
\omega_0=(1-p_1^2)dq_1\wedge dp_2+ p_1p_2(dq_1\wedge dp_1) +
(1+p_2^2)dq_2\wedge dp_1.
$$
with $q_1=t$ and $q_2=x$. A direct computation gives
$$
d\omega_0=3(p_1dp_2 - p_2dp_1)\wedge \Omega,
$$
and then the  Born - Infeld equation is not of divergent type.
\end{Ex}

\begin{Ex}[Tricomi equation]
The Tricomi equation is
$$
v_{xx} xv_{yy}+\alpha v_x + \beta v_y + \gamma(x,y).
$$
The corresponding primitive form is
$$
\omega_0=(\alpha p _1 + \beta p_2 +\gamma(q))dq_1\wedge dq_2+
dq_1\wedge dp_2-q_2dq_2\wedge dp_1,
$$
with $x=q_1$ and $y=q_2$. Since
$$
d\omega_0=(-\alpha dq_2 + \beta dq_1)\wedge \Omega,
$$
we conclude that the Tricomi equation is of divergent type.
\end{Ex}

\begin{Lem}
A Monge-Amp\`ere equation $\Delta_{\omega}=0$ is of divergent type
if and only if it exists a function $\mu$ on $M$ such that the
form $\omega + \mu \Omega$ is closed.
\end{Lem}

\begin{proof}
Since the exterior product by $\Omega$ is an isomorphism from
$\Omega^1(M)$ to $\Omega^3(M)$, for any $2$-form $\omega$, there
exists a $1$-form $\alpha_\omega$ such that
$$
d\omega=\alpha_\omega\wedge\Omega.
$$
Since $\bot(\alpha_\omega\wedge\Omega)=\alpha_\omega$ we deduce
that $\mathcal{E}(\omega)=0$ if and only if $d\alpha_\omega=0$,
that is $d(\omega+\mu\Omega)=0$ with $d\mu=-\alpha_\omega$.
\end{proof}

Hence, if $\Delta_\omega=0$ is of divergent type, one can choose
$\omega$ being closed.   The point is that it is not primitive in
general .

\subsection{Hitchin pairs}

Let us denote by $T$ the tangent bundle of $M$ and by $T^*$ its
cotangent bundle. The natural indefinite interior product on
$T\oplus T^*$ is
$$
(X+\xi,Y+\eta)=\frac{1}{2}(\xi(Y)+\eta(X)),
$$
and the Courant bracket on sections of $T\oplus T^*$ is
$$
[X+\xi,Y+\eta]=[X,Y]+L_X\eta-L_Y\xi -\frac{1}{2}d(\iota_X\eta-
\iota_Y\xi).
$$
\begin{Def}[Hitchin \cite{H1}]
An almost generalized complex structure is a bundle map
$\mathbb{J}: T\oplus T^*\rightarrow T\oplus T^*$ satisfying
$$
\mathbb{J}^2=-1,
$$
and
$$
(\mathbb{J}\cdot,\cdot)=-(\cdot,\mathbb{J}\cdot).
$$
Such an almost generalized complex structure is  said to be
integrable if the spaces of sections of its two eigenspaces are
closed under the Courant bracket.
\end{Def}

The standard examples are
$$
\mathbb{J}_1=\begin{pmatrix} J&0\\0&-J^*\end{pmatrix}
$$
and
$$
\mathbb{J}_2=\begin{pmatrix} 0&\Omega^{-1}\\ -\Omega &
0\end{pmatrix} $$ with $J$ a  complex structure and $\Omega$ a
symplectic form.

\begin{Lem}[Crainic \cite{Cr}]
Let $\Omega$ be a symplectic form and $\omega$ any $2$-form.
Define the tensor $A$ by $\omega=\Omega(A\cdot,\cdot)$ and the
form  $\tilde{\omega}$ by
$\tilde{\omega}=-\Omega(1+A^2\cdot,\cdot)$.

The almost generalized complex structure
\begin{equation}{\label{C}}
\mathbb{J}=\begin{pmatrix} A& \Omega^{-1}\\ \tilde{\omega} & -A^*
\end{pmatrix}
\end{equation}
is integrable if and only if $\omega$ is closed. Such a pair
$(\omega,\Omega)$ with $d\omega=0$ is called a Hitchin pair
\end{Lem}

We get then immediatly the following:

\begin{Prop}
To any $2$-dimensional symplectic Monge-Amp\`ere equation of
divergent type $\Delta_\omega=0$ corresponds a Hitchin pair
$(\omega,\Omega)$ and therefore a $4$-dimensional generalized
complex structure.
\end{Prop}

\begin{Rem}
Let $L^2\subset M^4$ be a $2$-dimensional submanifold. Let
$T_L\subset T$ be its tangent bundle and $T_L^0\subset T^*$ its
annihilator. $L$ is a generalized complex submanifold (according
to the terminology of \cite{Gu1}) or a generalized lagrangian
submanifold (according to the terminology of \cite{BB}) if
$T_L\oplus T^0_L$ is closed under $\mathbb{J}$. When $\mathbb{J}$
is defined by $\eqref{C}$, this is equivalent to saying that $L$
is lagrangian with respect to $\Omega$ and closed under $A$, that
is, $L$ is a generalized solution of $\Delta_\omega=0$.
\end{Rem}

\subsection{Conservation laws and Generating functions}

The notion of conservation laws is a natural generalization to
partial differential equations of the notion of  first integrals
(see \cite{KLR} for more details).

A $1$-form $\alpha$ is a conservation law for the equation
$\Delta_\omega=0$ if the restriction of $\alpha$ to  any
generalized solution is closed. Note that conservations laws are
actually well defined up closed forms.

\begin{Ex}
Let us consider the Laplace equation and the complex structure $J$
associated with. The $2$-form  $d\alpha$ vanish on any complex
curve if and only if $[d\alpha]_{1,1}=0$, that is
$$
\overline{\partial}\alpha_{1,0} + \partial \alpha_{0,1}=0
$$
or equivalently
$$
\overline{\partial}\alpha_{1,0}= \overline{\partial}\partial \psi
$$
for some real function $\psi$. (Here $\overline{\partial}$ is the
usual Dolbeault operator defined by the integrable complex
structure $J$.) We deduce that $\alpha -d\psi = \beta_{1,0}+
\beta_{0,1}$ with $\beta_{1,0}=\alpha_{1,0}-\partial\psi$ is a
holomorphic $(1,0)$-form.

Hence, the conservation laws of the $2D$-Laplace equation are (up
exact forms)  real part of $(1,0)$-holomorphic forms.
\end{Ex}

According to the Hodge-Lepage-Lychagin theorem, $\alpha$ is a
conservation law if and only if there exist two functions $f$ and
$g$ such that $d\alpha=f\omega+ g\Omega$. The function $f$ is
called a generating function of the Monge-Amp\`ere equation
$\Delta_\omega=0$. By analogy with the Laplace equation, we will
say that the function $g$  is the conjugate function to the
generating function $f$. We show in \cite{Ba3} that these
generating functions can be understood as "generalized harmonic
functions" for $2D$ Monge-Amp\`ere equations of divergent type.

The tensor $\mathbb{J}$ lives in $so(n,n)$, which can be
identified with the space of $2$-forms on $T\oplus T^*$, using the
inner product.  Moreover the space of forms, is isomorphic as a
Clifford algebra to the space of endomorphisms on $T$, and
therefore $\mathbb{J}$ acts on the tangent bundle. One obtains
then a differential operator $\bar{\partial}$, we can write as
follows:
$$
\bar{\partial} = d + i \mathbb{J}\circ d\circ \mathbb{J}.
$$
This operator is the exact analog to the Dolbeault operator in
complex geometry. M. Gualtieri actually proves in \cite{Gu1} that
an almost complex structure is integrable if and only if
$\bar{\partial}^2=0$.

Using this Gualtieri operator, we get the following
characterization for generating functions of symplectic $2D$ MAE
of divergent type:

\begin{Theo}
A function $f$ is a generating function of  the symplectic $2D$
MAE of divergent type $\Delta_\omega=0$ if and only if
$$
\partial_\omega\overline{\partial_\omega} f =0.
$$
\end{Theo}

\begin{Rem}[Reduction of the special lagrangian equation]

Only few explicit examples of special lagrangian submanifolds of
$\mathbb{C}^n$ are known and almost none in compact Calabi-Yau
manifolds (see papers of R. Bryant and D. Joyce on the subject)
and any explicit new example would be of great interest.

Historical first examples were constructed by R. Harvey and B.
Lawson in \cite{HL}. They considered invariant solutions with
respect to $SO(n)$ or a maximal torus. I propose to keep this
approach, considering action groups on the generalized complex
manifold $\mathbb{C}^n$. We know after \cite{BCG} that the reduce
space can admit also a generalized complex structure and when it
is $4$ dimensional, its generalized lagrangian submanifolds are
solutions of  $2D$-Monge-Amp\`ere equations of divergent type. The
idea would be the to construct global solutions as fibrations over
these Monge-Amp\`ere solutions.

If this approach is efficient, it should be generalized to other
classical examples of Calabi-Yau manifolds as $T^*S^n$ or degree
$n+2$ hypersurfaces in $P^{n+1}(\mathbb{C})$. This would be a
great motivations to investigate the geometry of
$2D$-Monge-Amp\`ere equations on $T^*S^2$ or $\mathbb{C}P^2$.
\end{Rem}

\end{document}